\documentclass[12pt]{article}

\usepackage{amsthm, amsmath}
\usepackage{amsfonts, amssymb}
\usepackage{mathrsfs, latexsym}

\newtheorem{thm}{Theorem}[section]
\newtheorem{prop}[thm]{Proposition}
\theoremstyle{definition}
\newtheorem{defn}{Definition}[section]
\theoremstyle{remark}

\theoremstyle{plain}

\newtheorem{cor}[thm]{Corollary}

\begin{document}
\title{Dimension of RC-lattices}
\maketitle
\begin{center}
\author{Dr. A. N. Bhavale}\\
{Head, Department of Mathematics,\\
PES Modern College of Arts, Science and Commerce (Autonomous),\\
Shivajinagar, Pune - 411005, M.S., India.\\
{\it email}: hodmaths@moderncollegepune.edu.in}
\end{center}
\begin{abstract}
In $1941$ Dushnik and Miller introduced the concept of dimension of a poset. 
In $2020$ Bhavale and Waphare introduced the concept of an RC-lattice as a lattice in which all the reducible elements are lying on a chain. 
In this paper, we introduce the concept of a complete fundamental basic block and prove that its dimension is at the most three. 
Consequently, we prove that the dimension of an RC-lattice on $n$ elements is at the most three. 
Further, we prove that an RC-lattice is non-planar if and only if its dimension is three.
\end{abstract}
\noindent
{\small {Keywords: }{Dimension, Poset, Lattice, Linear extension.}}\\
{\small {MSC classification 2020:} 06A05, 06A06, 06A07.}
\section{Introduction}
In $1930$ Szpilrajn \cite{SE} proved that any order relation on a set $X$ can be extended to a linear order on $X$. It also follows that any order relation is the intersection of its linear extensions. In $1941$ Dushnik and Miller \cite{DM} introduced the concept of the (order) dimension of poset $P$ as the minimum number of linear extensions of $P$ whose intersection is exactly $P$. It is known that the problem of finding dimension of posets is NP-complete.
In $1955$ Hiraguchi \cite{H} proved that the dimension of a poset does not exceed its width.
Hiraguchi \cite{H} also showed that the dimension of the poset $P$ on $n$ elements is $\leq \frac{n}{2}$. Komm \cite{K} showed that dimension of the poset consisting of all subsets of an $n$ element set ordered by inclusion is $n$.   For more results on dimension of posets, see  \cite{H1}, \cite{KR2}, \cite{TM}, and \cite{DK}. In this paper, we will restrict ourselves to the finite posets and lattices. 

In $2002$ Thakare, Pawar, and Waphare \cite{TPW} introduced the concept of an adjunct sum of two lattices.
In $2020$ Bhavale and Waphare \cite{BW} introduced the concepts of a nullity of a poset, a basic block associated to a poset, and a fundamental basic block associated to a dismantlable lattice.
Recently, Bhavale \cite{AB2} introduced the concept of a complete fundamental basic block on $r$ comparable reducible elements. 

In Section \ref{sec2}, we prove that the dimension of a complete fundamental basic block on $r$ comparable reducible elements is at the most three. Consequently, in Section \ref{sec4}, we prove that the dimension of an RC-lattice is at the most three. 
In Section \ref{sec3}, we obtain the dimension of an adjunct sum of two lattices in terms of the dimensions of the individual lattices. 
Moreover, we obtain a bound on the dimension of a dismantlable lattice in terms of its nullity.
In Section \ref{sec4}, we prove that the dimension of an RC-lattice is the same as the dimension of the basic block (and also the fundamental basic block) associated to that lattice. We also obtain a characterization, namely, an RC-lattice is of dimension three if and only if it is non planar.

A binary relation $\leq$ on a set $P$ is called a partial order if it is reflexive, antisymmetric and transitive on it. A partially ordered set (or \emph{poset}) consists of some ground set $P$ and a partial order $\leq$ defined on $P$. 
An element $c\in P$ is {\it a lower bound} ({\it an upper bound}) of $a,b\in P$ if $c\leq a, c\leq b (a\leq c, b\leq c)$. The {\it meet} of $a,b \in P$, is defined as the greatest lower bound of $a$ and $b$. The {\it join} of $a,b\in P$, is defined as the least upper bound of $a$ and $b$. A poset $L$ is a {\it lattice} if meet and join of $a$ and $b$ exist in $L$, $\forall ~ a,b\in L$.
A pair of elements $x, y \in P$ are \emph{comparable} if either $x \leq y$ or $y \leq x$. 
A pair of elements not comparable in $P$ are called \emph{incomparable}. If $x$ is incomparable to $y$ in $P$, we denote it by $x \parallel y$. 
A partial order $\leq$ on $P$ is a \emph{total order or linear order or chain} if for all $a, b \in P$, either $a \leq b$ or $b \leq a$. 
A partial order $\leq$ on $P$ is an \emph{antichain} if for all $a, b \in P$, $a \parallel b$. 
For $a \leq b$ in $P$, the interval $(a, b) = \{x \in P ~|~ a < x < b \}$.
For $a \in P$, let $(a] = \{x \in P ~|~ x \leq a\}$, and $D(a) = \{x \in P | x < a\}$.
Note that $D(a) \cup \{a\} = (a]$. Dually one can define $[a)$ and $U(a)$. Let $I(a) = \{b \in P ~|~ a \parallel b \}$.
An element $b$ in $P$ {\it covers} an element $a$ in $P$ if $a<b$, and there is no element $c$ in $P$ such that $a<c<b$. Denote this fact by $a\prec b$, and say that the pair $<a, b>$ is a $covering$ or an $edge$. If $a \prec b$ then $a$ is called a {\it lower cover} of $b$, and $b$ is called an {\it upper cover} of $a$. An element of a poset $P$ is called {\it doubly irreducible} if it has at most one lower cover and at most one upper cover in $P$. Let $Irr(P)$ denote the set of all doubly irreducible elements in the poset $P$.  Let $Red(P) = P \setminus Irr(P)$. The set of all coverings in $P$ is denoted by $E(P)$. The graph on a poset $P$ with edges as covering relations is called the {\it cover graph} and is denoted by $C(P)$. 
The {\it nullity of a graph} $G$ is given by $m-n+c$, where $m$ is the number of edges in $G$, $n$ is the number of vertices in $G$, and $c$ is the number of connected components of $G$. Bhavale and Waphare \cite{BW} defined {\it nullity of a poset} $P$, denoted by $\eta(P)$, to be the nullity of its cover graph $C(P)$. 
For an integer $n \geq 2$, a {\it crown} is a partially ordered set on $2n$ elements say, $x_1, y_1, x_2, y_2, \ldots , x_n, y_n$ in which $x_i\leq y_i$, 
$y_i\geq x_{i+1}$, for $i=1,2,\ldots,n-1$, and $x_1\leq y_n$ are the only comparability relations.

A total order is a \emph{linear extension} of a partial order $\leq$ if $\leq$ is a subset of it. 
A \emph{realizer} of a poset $P$ is a collection of linear extensions $R = \{R_1, R_2, \ldots , R_t\}$ of a partial order $\leq$ whose intersection is $\leq$. 
That is, for any incomparable pair $x, y \in P$, there are $R_i, R_j \in R$ with $(x, y) \in R_i$ and $(y, x) \in R_j$, i.e., $x \leq y$ in $R_i$ and $y \leq x$ in $R_j$. 
The \emph{dimension} of $P$ is defined as the size of a smallest realizer of $P$, and is denoted by $Dim(P)$. 
It is clear that the dimension of a chain is one, and that of an antichain (on at least two elements) is two. 
The dimension of crown on $2n$ elements is $3$ (see \cite{KR2}). 
For a linear extension $E$ of a poset $P$, and for a subposet $Q$ of $P$, let $E(Q)$ denote a part of $E$ containing a linear extension of $Q$, called a {\it partial linear extension} of $P$ (see \cite{DK}).

If $M$ and $N$ are two disjoint posets, the {\it{direct sum}} (see \cite{RS}), denoted by $M \oplus N$, is defined by taking the following order relation on $M\cup N:x\leq y$ if and only if $x,y\in M$ and $x\leq y$ in $M$, or $x,y\in N$ and $x\leq y$ in $N$, or $x\in M,y\in N$. 
\begin{prop} \cite{S} \label{dsp}
If $Q$ is a subposet of poset $P$ then $Dim(Q) \leq Dim(P)$.
\end{prop}
\begin{defn}\cite{R}
\textnormal{A finite lattice of order $n$ is called {\it{dismantlable}} if there exists a chain $L_{1} \subset L_{2} \subset \ldots\subset L_{n}(=L)$ of sublattices of $L$ such that $|L_{i}| = i$, for all $i$}.
\end{defn}
\begin{thm} \cite{KR1} \label{crown}
A finite lattice is dismantlable if and only if it contains no crown.
\end{thm}
The concept of {\it adjunct operation of lattices}, is introduced by Thakare, Pawar, and Waphare \cite{TPW}. 
Suppose $L_1$ and $L_2$ are two disjoint lattices and $(a, b)$ is a pair of elements in $L_1$ such that $a<b$ and $a\not\prec b$. Define the partial order $\leq$ on $L = L_1 \cup L_2$ with respect to the pair $(a,b)$ as follows: $x \leq y$ in L if $x,y \in L_1$ and $x \leq y$ in $L_1$, or $x,y \in L_2$ and $x \leq y$ in $L_2$, or $x \in L_1,$ $ y \in L_2$ and $x \leq a$ in $L_1$, or $x \in L_2,$ $ y \in L_1$ and $b \leq y$ in $L_1$. 

It is easy to see that $L$ is a lattice containing $L_1$ and $L_2$ as sublattices. The procedure for obtaining $L$ in this way is called {\it{adjunct operation (or adjunct sum)}} of $L_1$ with $L_2$. We call the pair $(a,b)$ as an {\it{adjunct pair}}, and $L$ as an {\it{adjunct}} of $L_1$ with $L_2$ with respect to the adjunct pair $\alpha = (a,b)$, and write $L = L_1 ]^b_a L_2$ or $L = L_1 ]_{\alpha} L_2$. A diagram of $L$ is obtained by placing a diagram of $L_1$ and a diagram of $L_2$ side by side in such a way that the largest element $1$ of $L_2$ is at a lower position than $b$ and the least element $0$ of $L_2$ is at a higher position than $a$, and then by adding the coverings $<1,b>$ and $<a,0>$. This clearly gives $|E(L)|=|E(L_1)|+|E(L_2)|+2$. 
A lattice $L$ is called an {\it{adjunct of lattices}} $L_0, L_1, \ldots, L_k$ if it is of the form $L = L_0 ]^{b_1}_{a_1} L_1 \cdots ]^{b_{k}}_{a_{k}} L_k$. 
\begin{thm}\cite{TPW}\label{dac} A finite lattice is dismantlable if and only if it is an adjunct of chains.
\end{thm}
\begin{thm} \cite{BW} \label{r-1}
A dismantlable lattice $L$ containing $n$ elements is of nullity $l$ if and only if $L$ is an adjunct of $l+1$ chains.
\end{thm}
Thakare, Pawar, and Waphare \cite{TPW} defined a \textit{block} as a finite lattice in which the largest element is join-reducible and the least element is meet-reducible.
\begin{defn}\cite{BW}\label{basicblock}
A poset $P$ is a {\it{basic block}} if it is one element or $Irr(P) = \emptyset$ or removal of any doubly irreducible element reduces nullity by one.
\end{defn}
\begin{defn}\cite{BW}\label{bbas}
A block $B$ is a {\it basic block associated to a poset} $P$ if $B$ is obtained from the basic retract associated to $P$ by successive removal of all the pendant vertices.
\end{defn}
\begin{thm} \cite{AB} \label{chbb} 
Let $B$ be an RC-lattice on $n$ elements. Then $B$ is a basic block having nullity $k$ if and only if 
$B = C_0 ]^{b_1}_{a_1} C_1 ]^{b_2}_{a_2} C_2 \cdots ]^{b_{k}}_{a_{k}} C_{k}$ with $a_i, b_i \in C_0$, satisfying 
(i) $|C_i| = 1$, for all $i, \; 1 \leq i \leq k$, (ii) $n = |C_0| + k$,  (iii) $Irr(B) = k+m$, where $m$ is the number of distinct adjunct pairs $(a_i, b_i)$ such that the interval $(a_i, b_i) \subseteq Irr(B)$, and (iv) $|C_0| = |Red(B)| + m$.
Further, if $|Red(B)|=r$ then $n=r+m+k$.
\end{thm}
From Theorem \ref{chbb}, it is clear that if $B$ is the basic block of nullity $k$ then 
$B = C ]_{a_1}^{b_1} \{c_1\} \cdots ]_{a_k}^{b_k} \{c_k\}$, where $C$ is a maximal chain containing all the reducible elements.

\begin{defn}\cite{BW}\label{FBB} 
A dismantlable lattice $B$ is said to be a {\it{fundamental basic block}} if it is a basic block and all the adjunct pairs in an adjunct representation of $B$ are distinct.
\end{defn}
\begin{defn}\cite{BW}\label{fbbas}
Let $L$ be a dismantlable lattice. Let $B$ be a basic block associated to $L$. If $B$ itself is a fundamental basic block then we say that $B$ is a fundamental basic block associated to $L$; otherwise, let $(a,b)$ be an adjunct pair in an adjunct representation of $B$. If the interval $(a,b) \subseteq Irr(B)$ then remove all but two non-trivial ears associated to $(a,b)$ in $B$; otherwise, remove all but one non-trivial ear associated to $(a,b)$ in $B$. Perform the operation of removal of non-trivial ears associated to $(a,b)$ for each adjunct pair $(a,b)$ in an adjunct representation of $B$. The resultant sublattice of $B$ is called a {\it fundamental basic block associated to} $L$.
\end{defn}
\begin{thm}\cite{AB} \label{redb}
Let $B$ be a basic block associated to a poset $P$. Then $Red(B) = Red(P)$ and $\eta(B) = \eta(P)$.
\end{thm}
Using Definition \ref{fbbas} and by Theorem \ref{redb}, we have the following result.
\begin{cor} \label{redfb}
Let $F$ and $B$ be the fundamental basic block and the basic block associated to a lattice $L$ respectively. 
Then $Red(F) = Red(B)= Red(L)$ and $\eta(F) \leq \eta(B)$.
\end{cor}

For the other definitions, notation, and terminology, see \cite{T}, \cite{GG}, \cite{DW}, and \cite{BW}.

\section{Complete fundamental basic block} \label{sec2}
Recently, Bhavale \cite{AB2} introduced the concept of a complete fundamental basic block on $r$ reducible elements which are all comparable, denoted by $CF(r)$, and it is defined as the fundamental basic block on $r$ reducible elements, having nullity $\binom{r}{2}$.  

\subsection{The lattice $L_r$} \label{subsec1}
In Subsection \ref{subsec2} we find the dimension of $CF(r)$. For that purpose, here we find that the recursive definition of $CF(r)$ is more helpful.
Hence we define $CF(r)$ recursively in the following, and we use the notation $L_r$ instead of $CF(r)$ hereafter. 
So suppose $L_1$ consists of a single element, say $a_1$. Now in the following, we give a recursive definition of $L_r$ for $r \geq 2$.
\begin{defn} \label{lr}
For $r \geq 2$, define $L_r = (L_{r-1} \oplus \{x_{r-1}\} \oplus \{b_{\binom{r}{2}}\}) ]_{\alpha_{{\binom{r}{2}} - (r-2)}} \\
\{c_{{\binom{r}{2}} - (r-2)}\} ]_{\alpha_{{\binom{r}{2}} - (r-2) + 1}} \{c_{{\binom{r}{2}} - (r-2) + 1}\} \cdots  ]_{\alpha_{\binom{r}{2}}} \{c_{\binom{r}{2}}\}$
where $\alpha_i = (a_i, b_i)$ for all $i, \; 1 \leq i \leq {\binom{r}{2}}$, 
$a_{{\binom{r-1}{2}} + 1} = a_{{\binom{r}{2}} - (r-2)} \prec x_1 \prec a_{{\binom{r}{2}} - (r-2) +1} \prec x_2 \prec \cdots \prec a_{{\binom{r}{2}} - 1} \prec x_{r-2} \prec a_{\binom{r}{2}}$, and $b_{{\binom{r}{2}} - (r-2)} = b_{{\binom{r}{2}} - (r-2) +1} = \cdots = b_{\binom{r}{2}}$.
\end{defn}
Note that $L_2$ is nothing but a diamond $M_2 = \{0=a_1, x_1, c_1, b_1=1\}$ with $x_1 \parallel c_1$, and the following figure (see Figure I) shows $L_r$ for $3 \leq r \leq 6$.

\begin{center}
%TeXCAD Picture [Nullity-3-Fig4-(3 to 6).pic]. Options:
%\grade{\on}
%\emlines{\off}
%\epic{\off}
%\beziermacro{\on}
%\reduce{\on}
%\snapping{\off}
%\pvinsert{% Your \input, \def, etc. here}
%\quality{8.000}
%\graddiff{0.005}
%\snapasp{1}
%\zoom{9.5137}
\unitlength 1.5mm % = 4.268pt
\linethickness{0.4pt}
\ifx\plotpoint\undefined\newsavebox{\plotpoint}\fi % GNUPLOT compatibility
\begin{picture}(80.726,65.695)(0,0)
\put(1.997,4.94){\circle{1.682}}
\put(1.997,10.932){\circle{1.682}}
\put(1.997,16.923){\circle{1.682}}
\put(1.997,22.914){\circle{1.682}}
\put(1.997,28.906){\circle{1.682}}
\put(7.988,10.932){\circle{1.682}}
\put(7.988,16.923){\circle{1.682}}
\put(7.988,22.914){\circle{1.682}}
\put(13.98,10.932){\circle{1.682}}
\put(13.98,22.914){\circle{1.682}}
\put(13.98,34.897){\circle{1.682}}
\put(19.971,4.94){\circle{1.682}}
\put(19.971,10.932){\circle{1.682}}
\put(19.971,16.923){\circle{1.682}}
\put(19.971,22.914){\circle{1.682}}
\put(19.971,28.906){\circle{1.682}}
\put(19.971,34.897){\circle{1.682}}
\put(19.971,40.888){\circle{1.682}}
\put(25.962,16.923){\circle{1.682}}
\put(25.962,22.914){\circle{1.682}}
\put(25.962,28.906){\circle{1.682}}
\put(31.954,22.914){\circle{1.682}}
\put(31.954,28.906){\circle{1.682}}
\put(31.954,34.897){\circle{1.682}}
\put(37.945,10.932){\circle{1.682}}
\put(37.945,22.914){\circle{1.682}}
\put(37.945,34.897){\circle{1.682}}
\put(37.945,46.88){\circle{1.682}}
\put(43.936,4.94){\circle{1.682}}
\put(43.936,10.932){\circle{1.682}}
\put(43.936,16.923){\circle{1.682}}
\put(43.936,22.914){\circle{1.682}}
\put(43.936,28.906){\circle{1.682}}
\put(43.936,34.897){\circle{1.682}}
\put(43.936,40.888){\circle{1.682}}
\put(43.936,46.88){\circle{1.682}}
\put(43.936,52.871){\circle{1.682}}
\put(49.928,16.923){\circle{1.682}}
\put(49.928,28.906){\circle{1.682}}
\put(49.928,40.888){\circle{1.682}}
\put(61.91,16.923){\circle{1.682}}
\put(61.91,22.914){\circle{1.682}}
\put(61.91,28.906){\circle{1.682}}
\put(61.91,34.897){\circle{1.682}}
\put(61.91,40.888){\circle{1.682}}
\put(61.91,46.88){\circle{1.682}}
\put(61.91,52.871){\circle{1.682}}
\put(67.902,4.94){\circle{1.682}}
\put(67.902,10.932){\circle{1.682}}
\put(67.902,16.923){\circle{1.682}}
\put(67.902,22.914){\circle{1.682}}
\put(67.902,28.906){\circle{1.682}}
\put(67.902,34.897){\circle{1.682}}
\put(67.902,40.888){\circle{1.682}}
\put(67.902,46.88){\circle{1.682}}
\put(67.902,52.871){\circle{1.682}}
\put(67.902,58.862){\circle{1.682}}
\put(67.902,64.854){\circle{1.682}}
\put(73.893,10.932){\circle{1.682}}
\put(73.893,22.914){\circle{1.682}}
\put(73.893,34.897){\circle{1.682}}
\put(73.893,46.88){\circle{1.682}}
\put(73.893,58.862){\circle{1.682}}
\put(79.885,28.906){\circle{1.682}}
\put(79.885,34.897){\circle{1.682}}
\put(79.885,40.888){\circle{1.682}}
\put(1.997,5.676){\line(0,1){4.625}}
\put(1.997,11.878){\line(0,1){4.31}}
\put(1.997,17.869){\line(0,1){4.204}}
\put(1.997,23.86){\line(0,1){4.31}}
\put(2.628,28.59){\line(1,-1){4.94}}
\put(2.628,17.554){\line(1,1){4.73}}
%\emline(2.523,16.397)(7.673,11.878)
\multiput(2.523,16.397)(.025624215,-.022486556){201}{\line(1,0){.025624215}}
%\end
%\emline(2.628,5.466)(7.463,10.406)
\multiput(2.628,5.466)(.022488988,.022977879){215}{\line(0,1){.022977879}}
%\end
%\emline(2.418,28.38)(7.463,17.764)
\multiput(2.418,28.38)(.022423803,-.047183419){225}{\line(0,-1){.047183419}}
%\end
%\emline(2.312,5.676)(7.463,16.292)
\multiput(2.312,5.676)(.022491123,.046359254){229}{\line(0,1){.046359254}}
%\end
\put(19.971,5.781){\line(0,1){4.52}}
\put(19.971,11.983){\line(0,1){4.099}}
\put(19.971,17.869){\line(0,1){4.204}}
\put(19.971,23.65){\line(0,1){4.52}}
\put(19.971,29.852){\line(0,1){4.099}}
\put(19.971,35.738){\line(0,1){4.204}}
%\emline(14.4,35.528)(19.341,40.573)
\multiput(14.4,35.528)(.022455655,.022933435){220}{\line(0,1){.022933435}}
%\end
%\emline(14.295,34.266)(19.235,29.431)
\multiput(14.295,34.266)(.022977879,-.022488988){215}{\line(1,0){.022977879}}
%\end
%\emline(14.505,23.44)(19.341,28.38)
\multiput(14.505,23.44)(.022488988,.022977879){215}{\line(0,1){.022977879}}
%\end
\put(14.4,22.284){\line(1,-1){4.94}}
%\emline(14.505,11.457)(19.446,16.397)
\multiput(14.505,11.457)(.022455655,.022455655){220}{\line(0,1){.022455655}}
%\end
%\emline(14.505,10.406)(19.235,5.676)
\multiput(14.505,10.406)(.022417161,-.022417161){211}{\line(0,-1){.022417161}}
%\end
%\emline(20.076,40.153)(25.647,23.86)
\multiput(20.076,40.153)(.022463361,-.065694735){248}{\line(0,-1){.065694735}}
%\end
%\emline(20.181,5.781)(25.752,22.179)
\multiput(20.181,5.781)(.022463361,.066118572){248}{\line(0,1){.066118572}}
%\end
\put(20.497,40.363){\line(1,-2){5.256}}
\put(20.602,17.554){\line(1,2){5.256}}
%\emline(20.392,28.38)(25.963,17.974)
\multiput(20.392,28.38)(.022463361,-.041959863){248}{\line(0,-1){.041959863}}
%\end
%\emline(20.497,5.571)(25.963,16.187)
\multiput(20.497,5.571)(.022493012,.04368835){243}{\line(0,1){.04368835}}
%\end
\put(44.042,5.886){\line(0,1){4.204}}
\put(44.042,11.983){\line(0,1){4.204}}
\put(44.042,17.974){\line(0,1){4.099}}
\put(44.042,23.86){\line(0,1){4.31}}
\put(44.042,29.957){\line(0,1){3.994}}
\put(44.042,35.843){\line(0,1){4.204}}
\put(44.042,41.729){\line(0,1){4.52}}
\put(44.042,47.931){\line(0,1){4.099}}
%\emline(38.366,47.616)(43.306,52.451)
\multiput(38.366,47.616)(.022977879,.022488988){215}{\line(1,0){.022977879}}
%\end
\put(38.366,46.354){\line(1,-1){4.94}}
%\emline(38.471,35.423)(43.516,40.363)
\multiput(38.471,35.423)(.022933435,.022455655){220}{\line(1,0){.022933435}}
%\end
%\emline(38.471,34.161)(43.306,29.431)
\multiput(38.471,34.161)(.02291532,-.022417161){211}{\line(1,0){.02291532}}
%\end
%\emline(38.471,23.44)(43.306,28.38)
\multiput(38.471,23.44)(.022488988,.022977879){215}{\line(0,1){.022977879}}
%\end
%\emline(38.471,22.389)(43.306,17.449)
\multiput(38.471,22.389)(.022488988,-.022977879){215}{\line(0,-1){.022977879}}
%\end
%\emline(38.471,11.562)(43.201,16.503)
\multiput(38.471,11.562)(.022417161,.023413479){211}{\line(0,1){.023413479}}
%\end
%\emline(32.374,35.528)(43.516,52.135)
\multiput(32.374,35.528)(.022463361,.033483123){496}{\line(0,1){.033483123}}
%\end
%\emline(32.374,34.161)(43.516,17.869)
\multiput(32.374,34.161)(.022463361,-.0328473675){496}{\line(0,-1){.0328473675}}
%\end
%\emline(32.374,23.44)(43.621,40.153)
\multiput(32.374,23.44)(.0224489793,.0333587636){501}{\line(0,1){.0333587636}}
%\end
%\emline(32.479,22.284)(43.621,5.676)
\multiput(32.479,22.284)(.022463361,-.033483123){496}{\line(0,-1){.033483123}}
%\end
%\emline(38.156,10.091)(43.201,5.571)
\multiput(38.156,10.091)(.025101272,-.022486556){201}{\line(1,0){.025101272}}
%\end
%\emline(44.357,5.571)(49.508,16.292)
\multiput(44.357,5.571)(.022491123,.046818257){229}{\line(0,1){.046818257}}
%\end
%\emline(44.462,28.38)(49.508,17.764)
\multiput(44.462,28.38)(.022423803,-.047183419){225}{\line(0,-1){.047183419}}
%\end
%\emline(44.462,17.554)(49.402,28.38)
\multiput(44.462,17.554)(.022455655,.049211329){220}{\line(0,1){.049211329}}
%\end
\put(44.357,40.258){\line(1,-2){5.256}}
%\emline(44.567,29.536)(49.508,40.153)
\multiput(44.567,29.536)(.022455655,.048255769){220}{\line(0,1){.048255769}}
%\end
%\emline(44.357,52.24)(49.613,41.624)
\multiput(44.357,52.24)(.022459738,-.045368672){234}{\line(0,-1){.045368672}}
%\end
%\emline(32.374,29.536)(43.832,52.135)
\multiput(32.374,29.536)(.0224650231,.0443117428){510}{\line(0,1){.0443117428}}
%\end
%\emline(32.269,28.275)(43.832,5.676)
\multiput(32.269,28.275)(.0224510162,-.0438815317){515}{\line(0,-1){.0438815317}}
%\end
\put(68.007,5.886){\line(0,1){4.31}}
\put(68.007,11.878){\line(0,1){4.31}}
\put(68.007,17.869){\line(0,1){4.099}}
\put(68.007,23.86){\line(0,1){4.415}}
\put(68.007,29.957){\line(0,1){4.099}}
\put(68.007,35.738){\line(0,1){4.31}}
\put(68.007,41.729){\line(0,1){4.415}}
\put(68,47.75){\line(0,1){4.5}}
\put(68,53.625){\line(0,1){4.375}}
\put(68,59.75){\line(0,1){4.375}}
%\emline(68.25,64.375)(73.25,59.625)
\multiput(68.25,64.375)(.023584799,-.022405559){212}{\line(1,0){.023584799}}
%\end
\put(68.375,53.375){\line(1,1){4.875}}
%\emline(68.25,52.5)(73.375,47.625)
\multiput(68.25,52.5)(.023617405,-.022465336){217}{\line(1,0){.023617405}}
%\end
\put(68.375,41.5){\line(1,1){4.875}}
%\emline(68.5,40.625)(73.5,35.75)
\multiput(68.5,40.625)(.02304137,-.022465336){217}{\line(1,0){.02304137}}
%\end
%\emline(68.5,29.375)(73.5,34.5)
\multiput(68.5,29.375)(.022421423,.022981959){223}{\line(0,1){.022981959}}
%\end
%\emline(68.25,28.5)(73.25,23.5)
\multiput(68.25,28.5)(.022421423,-.022421423){223}{\line(1,0){.022421423}}
%\end
%\emline(62,53.75)(67.125,64.625)
\multiput(62,53.75)(.022477968,.047697152){228}{\line(0,1){.047697152}}
%\end
\put(61.875,52.125){\line(1,-2){5.25}}
\put(61.875,41.75){\line(1,2){5.375}}
%\emline(61.75,40.125)(67.25,29.625)
\multiput(61.75,40.125)(.0224488779,-.0428569488){245}{\line(0,-1){.0428569488}}
%\end
%\emline(61.875,29.875)(67.25,40.5)
\multiput(61.875,29.875)(.022489438,.044455866){239}{\line(0,1){.044455866}}
%\end
\put(61.875,28.125){\line(1,-2){5.25}}
%\emline(68.375,17.5)(73.25,22.5)
\multiput(68.375,17.5)(.022465336,.02304137){217}{\line(0,1){.02304137}}
%\end
%\emline(62,17.875)(67.25,28.5)
\multiput(62,17.875)(.022435796,.045405777){234}{\line(0,1){.045405777}}
%\end
\put(61.875,16.125){\line(1,-2){5.25}}
%\emline(68.25,16.375)(73.625,12)
\multiput(68.25,16.375)(.027563978,-.022435796){195}{\line(1,0){.027563978}}
%\end
%\emline(68.375,5.5)(73.5,10.375)
\multiput(68.375,5.5)(.023617405,.022465336){217}{\line(1,0){.023617405}}
%\end
%\emline(61.875,22.125)(67.375,5.875)
\multiput(61.875,22.125)(.0224488779,-.0663262303){245}{\line(0,-1){.0663262303}}
%\end
%\emline(61.75,23.75)(67.5,40.375)
\multiput(61.75,23.75)(.0224608358,.0649411122){256}{\line(0,1){.0649411122}}
%\end
%\emline(61.875,34.125)(67.5,17.875)
\multiput(61.875,34.125)(.0224102571,-.0647407427){251}{\line(0,-1){.0647407427}}
%\end
%\emline(61.875,35.75)(67.5,52.125)
\multiput(61.875,35.75)(.0224102571,.0652387484){251}{\line(0,1){.0652387484}}
%\end
%\emline(61.875,46.25)(67.625,29.75)
\multiput(61.875,46.25)(.0224608358,-.0644528331){256}{\line(0,-1){.0644528331}}
%\end
%\emline(61.875,47.75)(67.625,64.25)
\multiput(61.875,47.75)(.0224608358,.0644528331){256}{\line(0,1){.0644528331}}
%\end
%\emline(68.125,5.625)(79.75,28.125)
\multiput(68.125,5.625)(.0224853914,.0435201124){517}{\line(0,1){.0435201124}}
%\end
%\emline(68.125,52.125)(79.75,29.75)
\multiput(68.125,52.125)(.0224853914,-.043278334){517}{\line(0,-1){.043278334}}
%\end
%\emline(68.125,17.875)(79.875,40.125)
\multiput(68.125,17.875)(.0224664375,.0425428284){523}{\line(0,1){.0425428284}}
%\end
%\emline(68.125,64)(79.875,41.75)
\multiput(68.125,64)(.0224664375,-.0425428284){523}{\line(0,-1){.0425428284}}
%\end
%\emline(68.125,64)(80,35.75)
\multiput(68.125,64)(.0224904285,-.0535035456){528}{\line(0,-1){.0535035456}}
%\end
%\emline(68,5.75)(80,34.375)
\multiput(68,5.75)(.0224718084,.0536046262){534}{\line(0,1){.0536046262}}
%\end
\put(1.997,.42){\makebox(0,0)[cc]{$L_3$}}
\put(19.971,.526){\makebox(0,0)[cc]{$L_4$}}
\put(44.042,.42){\makebox(0,0)[cc]{$L_5$}}
\put(68.007,.526){\makebox(0,0)[cc]{$L_6$}}
\end{picture}
\end{center}
\begin{center}
Figure I : Complete fundamental basic blocks
\end{center}

Observe that for each $i \geq 2$, $\eta(L_i) = \eta(L_{i-1}) + (i-1)$, since by Theorem \ref{r-1}, the nullity of a fundamental basic block is same as one less than the actual number of chains, i.e., the number of adjunct pairs, in its adjunct representation.
Therefore $\eta(L_r) = 1 + 2 + 3 + \cdots + (r-1) = {\binom{r}{2}}$.
Further, $a_1 \leq a_2 \leq \cdots \leq a_{\binom{r}{2}}$, and for $i < j$, $b_i \leq b_j$ whenever $a_i=a_j$.

For $r \geq 3$, let $A_1 = a_1 = a_2 = a_4 = a_7 = \cdots = a_{{\binom{r-2}{2}} + 1} = a_{{\binom{r-1}{2}} + 1}$,
for $2 \leq k \leq r-1$, let $A_k = a_{\binom{k+1}{2}} = a_{{\binom{k+1}{2}} + (k)} = a_{{\binom{k+1}{2}} + (k) + (k+1)} = \cdots = 
a_{{\binom{r-1}{2}} + k} = b_{{\binom{k}{2}} - (k-2)} = b_{{\binom{k}{2}} - (k-2) + 1} = \cdots = b_{{\binom{k}{2}} - 1} = b_{\binom{k}{2}}$, and 
let $A_r = b_{{\binom{r-1}{2}} + 1} = b_{{\binom{r-1}{2}} + 2} = \cdots = b_{{\binom{r-1}{2}} + (r-1)}$.

Note that ${\binom{r-1}{2}} + (r-1) = {\binom{r}{2}}$ and ${\binom{k+1}{2}} + (k) + (k+1) + \cdots + (k+(r-k-2)) = {\binom{r-1}{2}} + k$.
Also for $1 \leq k \leq r$, there are $(r-k) \; a_i's$ which are same as that of $A_k$, and there are $(k-1) \; b_i's$ which are same as that of $A_k$. 
For $r \geq 2$, let $C : A_1 \prec x_1 \prec A_2 \prec x_2 \prec A_3 \cdots  \prec x_{r-2} \prec A_{r-1} \prec x_{r-1} \prec A_r$.
Then $C$ is a maximal chain in $L_r$ containing all the $r$ reducible elements $A_1, A_2, \ldots , A_r$.
Thus it follows that $L_r$ is the fundamental basic block containing $r$ reducible elements which are all comparable.
Also the fundamental basic block containing $r$ comparable reducible elements with the largest nullity is $L_r$.
It can be observed that for all $1 \leq i \leq r$, $(A_i] \cong L_i$ and $[A_i) \cong L_{r-i+1}$.
Also $|L_r| = \frac{r^2+3r-2}{2}$, since $|L_r| = |C| + {\binom{r}{2}} = (2r-1) + {\binom{r}{2}}$.
Also for $\binom{r-1}{2} + 1 = \binom{r}{2} - (r-2) \leq i \leq \binom{r}{2}$, $a_i \prec c_i \prec b_i$.
Therefore $A_{i+1} \prec c_{\binom{k}{2}+1+i} \prec A_{k+1}$ for all $i$, $0 \leq i \leq k-1$.

\subsection{Dimension of $L_r$} \label{subsec2}
In order to obtain a bound for the dimension of an RC-lattice, firstly we obtain the dimension of $L_r$.
For that purpose, we take the help of the following result due to Baker et al. \cite{BFR}.
\begin{thm} \cite{BFR} \label{chpl}
The dimension of a lattice is at most two if and only if it is planar.
\end{thm}
Therefore by Theorem \ref{chpl}, it is clear that $Dim(L_1) = 1$, and for $r = 2, 3, 4, \; Dim(L_r) = 2$, since they are all planar lattices.
\begin{thm} \label{DLr3}
For $r \geq 5$, $Dim(L_r) = 3$.
\end{thm}
\begin{proof}
Let $r \geq 5$. Clearly $L_5$ is non planar (see Figure I). Therefore by Definition \ref{lr}, $L_r$ is non planar for all $r \geq 5$. 
Hence $Dim(L_r) \geq 3$.
Now for $r \geq 5$, let $R_1^r$ be the linear extension $A_1 \prec c_{{\binom{r-1}{2}} +1} \prec \cdots \prec c_4 \prec c_2 \prec c_1 \prec x_1 \prec
          A_2 \prec c_{{\binom{r-1}{2}} +2} \prec \cdots \prec c_8 \prec c_5 \prec c_3 \prec x_2 \prec
          A_3 \prec c_{{\binom{r-1}{2}} +3} \prec \cdots \prec c_{13} \prec c_9 \prec c_6 \prec x_3 \prec \cdots \prec
          A_i \prec c_{{\binom{r-1}{2}} +i} \prec \cdots \prec c_{{\binom{i+1}{2}} + (i) + (i+1)}  \prec c_{{\binom{i+1}{2}} + (i)}
                \prec c_{{\binom{i+1}{2}}} \prec x_i \prec  \cdots \prec A_{r-1} \prec 
		    c_{{\binom{r-1}{2}} + (r-1)} \prec x_{r-1} \prec A_r$,
$R_2^r$ be the linear extension $A_1  \prec x_1  \prec c_1  \prec
          A_2  \prec x_2  \prec c_3 \prec c_2 \prec
          A_3  \prec x_3  \prec c_6 \prec c_5 \prec c_4 \prec \cdots \prec
          A_{r-1}  \prec x_{r-1} \prec c_{{\binom{r}{2}}} \prec c_{{\binom{r}{2}} - 1} \prec c_{{\binom{r}{2}} - 2} \prec \cdots 
                        \prec c_{{\binom{r}{2}} - (r-2)} \prec A_r$, and					 							
$R_3^r$ be the linear extension $A_1  \prec c_{{\binom{r}{2}} - (r-2)} \prec x_1  \prec c_1  \prec
          A_2  \prec c_{{\binom{r}{2}} - (r-2) + 1}  \prec x_2  \prec c_3 \prec c_2 \prec
          A_3  \prec c_{{\binom{r}{2}} - (r-2) + 2} \prec x_3  \prec c_6 \prec c_5 \prec c_4 \prec \cdots \prec
          A_{r-2}  \prec c_{{\binom{r}{2}} - (r-2) + (r-3)} \prec x_{r-2} \prec c_{{\binom{r-1}{2}}} \prec c_{{\binom{r-1}{2}} - 1} \prec       
                        c_{{\binom{r-1}{2}} - 2} \prec \cdots  \prec c_{{\binom{r-1}{2}} - (r-3)} \prec 
          A_{r-1}  \prec c_{{\binom{r-1}{2}} + (r-1)}  \prec x_{r-1} \prec A_r$.
Clearly $R_1^r, R_2^r, R_3^r$ are linear extensions of $L_r$ for all $r \geq 5$.\\
Claim : For $r \geq 5$, $\{ R_1^r, R_2^r, R_3^r \}$ is a realizer of $L_r$. \\
Using the method of induction on $r$, we have for $r = 5$, 
$R_1^5 : A_1 \prec c_7 \prec c_4 \prec c_2 \prec c_1 \prec x_1 \prec A_2  \prec c_8 \prec c_5 \prec c_3 \prec x_2 \prec
          A_3 \prec c_9 \prec c_6 \prec x_3 \prec A_4 \prec c_{10} \prec x_4 \prec A_5$,
$R_2^5 : A_1 \prec x_1 \prec c_1 \prec A_2 \prec x_2 \prec c_3 \prec c_2  \prec
          A_3  \prec x_3  \prec c_6 \prec c_5 \prec c_4 \prec A_4  \prec x_4 \prec c_{10} \prec c_9 \prec c_8 \prec c_7  \prec A_5$, and 
$R_3^5 : A_1 \prec c_7 \prec x_1 \prec c_1 \prec A_2  \prec c_8 \prec x_2 \prec c_3 \prec c_2 \prec  
          A_3 \prec c_9 \prec x_3 \prec c_6 \prec c_5 \prec c_4 \prec A_4 \prec c_{10} \prec x_4 \prec A_5$.
It can be easily checked that $\{ R_1^5, R_2^5, R_3^5 \}$ is a realizer of $L_5$. 
Now assume that $\{ R_1^k, R_2^k, R_3^k \}$ is a realizer of $L_k$ for some $k \geq 5$. We prove that 
$\{ R_1^{k+1}, R_2^{k+1}, R_3^{k+1} \}$ is a realizer of $L_{k+1}$.
Let $A = \{x_k, A_{k+1}, c_{\binom{k+1}{2}}\}$ and $B = \{c_{{\binom{k}{2}} + 1}, c_{{\binom{k}{2}} + 2}, \ldots , c_{{\binom{k}{2}} + k}\}$.
Let $C_i = (A_{i+1}, A_{k+1}) \setminus \{c_{{\binom{k}{2}} + 1 + i} \}$ for all $i, ~ 0 \leq i \leq k-1$.
Let $D_i = L_{k+1} \setminus ((A_{i+1}] \cup U(A_{i+1}))$ for all $i, ~ 1 \leq i \leq k-1$.
Now for $k \geq 5$, observe that
\begin{enumerate}
\item $x \leq y$ in $L_{k+1}$ if and only if $x, y \in L_k$ with $x \leq y$, or $x \in L_k$ and $y \in A$, or $x \in B \cup \{x_k\}$ and $y = A_{k+1}$.
\item $x \parallel y$ in $L_{k+1}$ if and only if $x, y \in L_k$ with $x \parallel y$, or $x, y \in B$, or
$x = c_{{\binom{k}{2}} + 1 + i} \in B$ and $y \in C_i$ for all $i, ~ 0 \leq i \leq k-1$, or 
$x = c_{{\binom{k}{2}} + 1 + i} \in B$ and $y \in D_i$ for all $i, ~ 1 \leq i \leq k-1$. 
\end{enumerate}
Note that $A_{i+1} \prec c_{\binom{k}{2}+1+i} \prec A_{k+1}$ for all $i$, $0 \leq i \leq k-1$, and \\
$D_1 = \{ c_2, c_4, c_7, c_{11}, \ldots , c_{{\binom{k+1}{2}} - (k-1)} \}$, 
$D_2 = \{ c_4, c_5, c_7, c_8, c_{11}, c_{12}, \ldots , \\
c_{{\binom{k+1}{2}} - (k-1)}, c_{{\binom{k+1}{2}} - (k-1) + 1}\}$, 
$D_3 = \{ c_7, c_8, c_9, c_{11}, c_{12}, c_{13}, \ldots , c_{{\binom{k+1}{2}} - (k-1)}, \\
c_{{\binom{k+1}{2}} - (k-1) + 1}, c_{{\binom{k+1}{2}} - (k-1) + 2}\}$, $\ldots$ , 
$D_{k-1} = \{c_{{\binom{k+1}{2}} - (k-1)}, c_{{\binom{k+1}{2}} - (k-1) + 1}, \ldots , \\
c_{{\binom{k+1}{2}} - (k-1) + (k-2)}\}$.

It is clear from Definition \ref{lr} that the set $L_{k+1} \setminus (A \cup B) = L_k$.
Now the linear extensions $R_1^{k+1}, R_2^{k+1}, R_3^{k+1}$ restricted to $L_k$ are respectively the linear extensions $R_1^k, R_2^k, R_3^k$. Therefore by induction hypothesis, it is sufficient to take the care of all comparabilities and all incomparabilities of the elements in the sets $A$ and $B$ in $L_{k+1}$. 
By taking the care of all comparabilities means, $x \leq y$ in $L_{k+1}$ if and only if $x \leq y$ in each $R_i^{k+1}$ for $i=1,2,3$.
Similarly, by taking the care of all incomparabilities means, $x \parallel y$ in $L_{k+1}$ if and only if $x \leq y$ in $R_i^{k+1}$ and
$y \leq x$ in $R_i^{k+1}$ for $i \neq j$.

Now the care of both the observations mentioned above is taken by all the linear extensions $R_1^{k+1}$, $R_2^{k+1}$, and $R_3^{k+1}$.
Hence $\{R_1^{k+1}, R_2^{k+1}, R_3^{k+1}\}$ is a realizer of $L_{k+1}$. Thus by mathematical induction, $Dim(L_r) \leq 3$ for all $r \geq 5$.
Hence $Dim(L_r) = 3$ for all $r \geq 5$.
\end{proof}
As $Dim(L_1) = 1$ and $Dim(L_r) = 2$ for $r = 2, 3, 4$, by Theorem \ref{DLr3}, we have the following result.
\begin{thm} \label{dlr3}
For $r \geq 1$, $Dim(L_r) \leq 3$.
\end{thm}

\section{Dimension of lattices} \label{sec3}

The following result follows from Proposition \ref{dsp}, which gives a lower bound for the dimension of adjunct sum of two lattices.
\begin{prop} \label{dim}
Let $L_1, L_2$ be the lattices of dimension $m, n$ respectively. Let $L = L_1 ]^{b}_{a} L_2$. Then $Dim(L) \geq max \{m, n\}$.
\end{prop}
In the following result, we obtain an upper bound for the dimension of adjunct sum of two lattices.
\begin{prop} \label{dim1}
Let $L_1, L_2$ be the lattices of dimension $m, n$ respectively. Let $L = L_1 ]^{b}_{a} L_2$. Then $Dim(L) \leq max \{2m, m + n\}$, that is, $Dim(L) \leq m + max \{m, n\}$.
\end{prop}
\begin{proof}
Let $\{E_1, E_2, \ldots , E_m\}$ be a realizer of $L_1$. Let $\{F_1, F_2, \ldots , F_n\}$ be a realizer of $L_2$.
Then we have the following two cases.
Case 1: Suppose $m \geq n$.
Consider the linear extensions $E_i((a]) \oplus F_i \oplus E_i(U(a))$ for $1 \leq i \leq m$, where $F_i = F_1$ for $n+1 \leq i \leq m$,
and $E_i(D(b)) \oplus F_1 \oplus E_i([b))$ for $1 \leq i \leq m$.
As these $2m$ extensions forms a realizer of $L$, $Dim(L) \leq 2m$.
Case 2: Suppose $m < n$.
Consider the linear extensions $E_i((a]) \oplus F_i \oplus E_i(U(a))$ for $1 \leq i \leq n$, where $E_i = E_1$ for $m+1 \leq i \leq n$,
and $E_i(D(b)) \oplus F_1 \oplus E_i([b))$ for $1 \leq i \leq m$.
As these $n + m$ extensions forms a realizer of $L$, $Dim(L) \leq n + m$.
Thus $Dim(L) \leq max \{2m, m + n\}$.
\end{proof}
The more precise bound for the dimension of adjunct sum of two lattices is achieved in the following result.
\begin{thm} \label{dmn}
Let $L_1, L_2$ be the lattices of dimension $m, n$ respectively. Let $L = L_1 ]^{b}_{a} L_2$. Then 
$m \leq Dim(L) \leq m + 1$ for $m \geq n$, and $Dim(L) = n$ for $m < n$.
\end{thm}
\begin{proof}
Let $\{E=E_1, E_2, \ldots , E_m\}$ be a realizer of $L_1$. Let $\{F_1, F_2, \ldots , F_n\}$ be a realizer of $L_2$. 
Let $I(L_2) = \{ x \in L_1 ~|~ x \parallel y$ in $L, \; \forall \; y \in L_2 \}$. Define two extensions of $L$ as follows.
Let $M_1 = E((a] \cup I(L_2)) \oplus F_1 \oplus E([b))$, and let $M_2 = E((a]) \oplus F_2 \oplus E(I(L_2) \cup [b))$ (Take $F_2 = F_1$, if $n = 1$).
It is clear that $x \leq y$ in $L$ if and only if $x \leq y$ in $M_1$ and $x \leq y$ in $M_2$.
Also, $x \parallel y$ in $L$ if and only if $x, y \in L_1$ with $x \parallel y$, or $x, y \in L_2$ with $x \parallel y$, or 
$x \in L_1 \setminus ((a] \cap [b))$ and $y \in L_2$.
Now if $x \in L_1 \setminus ((a] \cap [b))$ and $y \in L_2$ then $x \leq y$ in $M_1$ and $y \leq x$ in $M_2$.
Therefore, if $m \geq n$ then the partial extensions $M_1, M_2, E_2, E_3, \ldots , E_m$ forms a realizer of $L$. Therefore $Dim(L) \leq m + 1$.
Also, if $m < n$ then the partial extensions $M_1, M_2, F_3, F_4, \ldots , F_n$ forms a realizer of $L$. Therefore $Dim(L) \leq n$.
Thus the proof follows from Proposition \ref{dim}.
\end{proof}
Using the above Theorem \ref{dmn}, we have the following result.
\begin{cor} \label{dmn1}
Let $L$ be a lattice with $|L| \geq 3$, and let $C$ be a chain. Let $L' = L ]_{a}^{b} C$. Then $Dim(L) = m$ implies that $m \leq Dim(L') \leq m+1$.
\end{cor}
Thus, the dimension of a lattice increases by at most one if a chain is added to it by means of an adjunct sum.
\begin{cor} \label{dmn2}
If $L = C_0 ]_{a_1}^{b_1} C_1 \cdots ]_{a_k}^{b_k} C_{k}$ where $C_i$ are chains for $0 \leq i \leq k$ then $Dim(L) \leq k+1$.
\end{cor}
Using Theorem \ref{r-1}, in the following result, we obtain a relation between the dimension and the nullity of a dismantlable lattice.
\begin{cor}
If $L$ is a dismantlable lattice of nullity $k$ then $Dim(L) \leq k+1$.
\end{cor}
It is clear that the dimension of adjunct sum of two chains is two.
In the following result, we prove that the dimension of a dismantlable lattice remains the same if a chain is added to it by means of an adjunct sum under the given restriction.
\begin{thm} \label{LC}
Let $L$ be a dismantlable lattice which is not a chain.
If $L' = L ]_a^b C$ where $C$ is a chain and $(a, b)$ is an adjunct pair of $L$ then $Dim(L') = Dim(L)$.
\end{thm}
\begin{proof}
Suppose $Dim(L) = n$. As $L$ is not a chain, $n \geq 2$. As $L$ is a sublattice of $L'$, by Proposition \ref{dsp}, $Dim(L') \geq n$.
Now let $\{R_1, R_2, \cdots , R_n\}$ be a realizer of $L$. 
Let $R'_1 = (R_1 \cap (a] ) \oplus C \oplus (R_1 \cap (U(a)))$, and for $2 \leq i \leq n$, let $R'_i = (R_i \cap (D(b)) ) \oplus C \oplus (R_i \cap [b))$.
Then $\{R'_1, R'_2, \cdots , R'_n\}$ is a realizer of $L'$. Hence $Dim(L') \leq n$. Thus $Dim(L') = n = Dim(L)$.
\end{proof}

\section{Dimension of RC-lattices} \label{sec4}
In this section, we prove that the dimension of an RC-lattice is same as the dimension of the basic block (and also the fundamental basic block) associated to that lattice.
Further, we prove that the dimension of an RC-lattice is at the most three.
By Theorem \ref{crown}, it is clear that an RC-lattice is dismantlable, as it does not contain a crown.
\begin{thm} \label{BL}
Let $L$ be an RC-lattice. If $B$ is the basic block associated to $L$ then $Dim(L) = Dim(B)$.
\end{thm}
\begin{proof}
As all the reducible elements in $L$ are comparable, $L$ is dismantlable lattice, as it does not contain a crown (see Theorem \ref{crown}).
By Theorem \ref{dac}, suppose $L = C_0 ]_{a_1}^{b_1} C_1 \dots ]_{a_k}^{b_k} C_{k}$, where $C_0$ is a maximal chain containing all the reducible elements, and $C_i$ is a chain for all $1 \leq i \leq k$. 
If $B$ is the basic block associated to $L$ then by Theorem \ref{chbb}, 
$B = C ]_{a_1}^{b_1} \{c_1\} \cdots ]_{a_k}^{b_k} \{c_k\}$, where $C$ is a maximal chain containing all the reducible elements, and $c_i \in C_i$ for all $1 \leq i \leq k$.
By Proposition \ref{dsp}, as $B$ is a sublattice of $L$, $Dim(B) \leq Dim(L)$.
Now we claim that $Dim(L) \leq Dim(B)$.
Let $Dim(B) = n$. Let $R = \{R_1, R_2, \ldots , R_n\}$ be a realizer of $B$.
For each $i,\; 1 \leq i \leq n$, let $R'_i$ be the linear extension obtained from $R_i$ as follows.
\begin{enumerate}
\item Replace each $c_i$ in $R_i$ by the chain $C_i$ for all $i, \; 1 \leq i \leq k$.
\item If $a < b$ are consecutive reducible elements of $L$, and $x \in C \cap Irr(B)$ with $a \prec x \prec b$ in $B$, then 
replace $x$ in $R_i$ by the chain $C_0 \cap (a, b)$ of $L$.
\item If $a < b$ are consecutive reducible elements of $L$, and $C \cap (a, b) = \emptyset$ in $B$, then 
between $a$ and $b$ in $R_i$, put the chain $C_0 \cap (a, b)$ of $L$.
\end{enumerate}
Then it follows that $\{R'_i ~|~ 1 \leq i \leq n\}$ is a realizer of $L$. Therefore $Dim(L) \leq Dim(B)$. 
Thus $Dim(L) = Dim(B)$.
\end{proof}
As $Dim(L_r) = 2$ for $r = 2, 3, 4$, by Lemma \ref{LC}, we have the following.
\begin{cor}
For $2 \leq r \leq 4$, if $L = L_r ]_a^b C$ where $C$ is a chain, and $(a, b)$ is an adjunct pair in $L_r$ then $Dim(L) = 2$.
\end{cor}
As $Dim(L_r) = 3$ for $r \geq 5$, by Lemma \ref{LC}, we have the following result.
\begin{cor}
For $r \geq 5$, if $L = L_r ]_a^b C$ where $C$ is a chain, and $(a, b)$ is an adjunct pair in $L_r$ then $Dim(L) = 3$.
\end{cor}
\begin{thm} \label{LBF}
If $F$ is the fundamental basic block associated to an RC-lattice $L$ then $Dim(F) = Dim(L)$.
\end{thm}
\begin{proof}
Let $B$ be the basic block associated to an RC-lattice $L$. Suppose $\eta(L) = k$. Then by Theorem \ref{redb}, $\eta(B) = k$. 
Therefore by Theorem \ref{chbb}, $B = C ]_{a_1}^{b_1} \{c_1\} ]_{a_2}^{b_2} \{c_2\} \cdots ]_{a_k}^{b_k} \{c_k\}$, where $C$ is a maximal chain containing all the reducible elements. Now $F$ is also the fundamental basic block associated to $B$. Therefore by Corollary \ref{redfb}, $\eta(F) \leq k$. Suppose $\eta(F) = l$ and $m = k - l$. Then 
by Theorem \ref{r-1}, and using Definition \ref{FBB} and Definition \ref{fbbas}, $F = C ]_{a_{i_1}}^{b_{i_1}} \{c_{i_1}\} ]_{a_{i_2}}^{b_{i_2}} \{c_{i_2}\} \cdots ]_{a_{i_l}}^{b_{i_l}} \{c_{i_l}\}$, where $i_j \in \{1, 2, \ldots , k\}$ for each $j, ~ 1 \leq j \leq l$, 
and $(a_{i_p}, b_{i_p}) \neq (a_{i_q}, b_{i_q})$ for $1 \leq p \neq q \leq l$.
As $B$ can be obtained from $F$ by taking adjunct of $F$ with $m$ chains (which consists of singletons), by Lemma \ref{LC}, $Dim(B) = Dim(F)$. Also by Theorem \ref{BL}, $Dim(L) = Dim(B)$. Hence $Dim(L) = Dim(F)$.
\end{proof}
In $1951$ Hiraguchi \cite{H} proved that removal of a chain decreases the dimension of a poset by at most two. In that regard we have the following.
\begin{prop}
If $L = C_0 ]_{a_1}^{b_1} C_1 \cdots ]_{a_k}^{b_k} C_{k}$ where $C_0$ is a maximal chain containing all the reducible elements, then 
$Dim(L) \leq Dim(L \setminus C_i) + 1, \; \forall \; i, \; 1 \leq i \leq k$.
\end{prop}
\begin{proof}
If $B$ is a basic block associated to $L$ then 
by Lemma \ref{BL}, $Dim(B) = Dim(L)$. By Theorem \ref{chbb}, suppose
$B = C ]_{a_1}^{b_1} \{c_1\} \dots ]_{a_k}^{b_k} \{c_k\}$ where $C$ is a maximal chain containing all the 
reducible elements, and $c_i \in C_i$ for $1 \leq i \leq k$.  
We know that removal of an element from a poset decreases its dimension by one (see \cite{H}).
Therefore $Dim(B) \leq Dim(B \setminus \{c_i\}) + 1, \; \forall \; i, \; 1 \leq i \leq k$.
For each $i, \; 1 \leq i \leq k$, if $B'_i$ is the basic block associated to $B \setminus \{c_i\}$ then it is also the basic block associated to $L \setminus C_i$. 
By Theorem \ref{BL}, $Dim(B'_i) = Dim(B \setminus \{c_i\}) = Dim(L \setminus C_i)$. Hence the proof.
\end{proof}
We now prove the main result in the following.
\begin{thm} \label{dl3}
If $L$ is an RC-lattice then $Dim(L) \leq 3$.
\end{thm}
\begin{proof}
Let $L$ be an RC-lattice. Let $B$ and $F$ be the basic block and the fundamental basic block associated to $L$ respectively. Then by Theorem \ref{BL} and by Theorem \ref{LBF}, $Dim(L) = Dim(B) = Dim(F)$. 
If $|Red(L)| = r$ then by Theorem \ref{redb} and Corollary \ref{redfb}, $|Red(B)| = |Red(F)| = r$.
Therefore $F$ is a sublattice of $L_r$, since $L_r$ is the fundamental basic block containing $r$ comparable reducible elements and having the largest nullity $\binom{r}{2}$.
Therefore by Proposition \ref{dsp}, $Dim(F) \leq Dim(L_r)$. Thus the proof follows from Theorem \ref{dlr3}.
\end{proof}
Finally, as a consequence of Theorem \ref{chpl} and Theorem \ref{dl3}, we have the following result.
\begin{thm}
Let $L$ be an RC-lattice. Then $Dim(L) = 3$ if and only if $L$ is non planar.
\end{thm}


\begin{thebibliography}{10}
\bibitem{SE} E. Szpilrajn, {\it Sur l'extension de l'ordre partiel}, Fundamenta Mathematicae, {\bf 16}(1930),386--389, doi:10.4064/fm-16-1-386-389.
\bibitem{DM} B. Dushnik and E. W. Miller, {\it Partially ordered sets}, American J. Math., \textbf{63}(1941), 600--610.
\bibitem{K} H. Komm, {\it On the dimension of partially ordered sets}, American J. Math., \textbf{20}(1948), 507-520.
\bibitem{H1} T. Hiraguchi, {\it On the dimension of partially ordered sets}, Sci. Rep. Kanazawa Univ., \textbf{1}(1951), 77--94.
\bibitem{H} T. Hiraguchi, {\it On the dimension of orders}, Sci. Rep. Kanazawa Univ., \textbf{4}(1955), 1--20.
\bibitem{BFR} K. A. Baker, P. C. Fishburn,  and F. S. Roberts, {\it Partial orders of dimension 2}, Networks, \textbf{2}(1971), 11--28.
\bibitem{R} I. Rival, {\it{Lattices with doubly irreducible elements}}, Canad. Math. Bull., \textbf{17}(1974), 91--95.
\bibitem{KR1} D. Kelly and I. Rival,  {\it Crown, fences and dismantlable lattices}, Canad. J. Math., \textbf{27} (1974), 636--665.
\bibitem{KR2} D. Kelly and I. Rival,  {\it Certain partially ordered sets of dimension three}, J. Comb. Theory (A), \textbf{18} (1975), 239--242.
\bibitem{TM} W. T. Trotter, Jr. and J. I. Moore, Jr.,  {\it The dimension of planar posets}, J. Comb. Theory (A), \textbf{22} (1977), 54--67.
\bibitem{DK} D. Kelly,  {\it On the dimension of partially ordered sets}, Discrete Mathematics, \textbf{35} (1981), 135--156.
\bibitem{RS} R. P. Stanley, {\it{Enumerative Combinatorics}}, Wadsworth and Brooks, California, 1986.
\bibitem{T} W. T. Trotter, Jr., {\it Combinatorics and partially ordered sets: Dimension theory}, Baltimore: Johns Hopkins University press, 1992.
\bibitem{GG} G. Gr\"{a}tzer, {\it{General Lattice Theory}}, Birkh\"{a}user Verlag, Second Ed., 1998.
\bibitem{TPW} N. K. Thakare, M. M. Pawar, and B. N. Waphare, {\it A structure theorem for dismantlable lattices and enumeration}, Periodica Mathematica Hungarica, \textbf{45}(2002), 147--160.
\bibitem{DW} D. B. West, {\it{Introduction to Graph Theory}}, Prentice Hall of India, New Delhi, 2002.
\bibitem{S}  B. S. W. Schr\"{o}der, Ordered Sets - An Introduction, {\it Birkh\"{a}user, Boston}, (2003).
\bibitem{BW} A. N. Bhavale and B. N. Waphare, {\it Basic retracts and counting of lattices}, Australas. J. Combin., {\bf 78(1)}(2020), 73--99.
\bibitem{AB} A. N. Bhavale, {\it Enumeration of lattices of nullity $k$ and containing $r$ comparable reducible elements}, Soft Computing, (2024) (Communicated).
https://doi.org/10.48550/arXiv.2502.06912 
\bibitem{AB2} A. N. Bhavale, {\it Equivalence of labeled graphs and lattices}, Australas. J. Combin., (2025) (Communicated). https://doi.org/10.48550/arXiv.2501.05064
\end{thebibliography}
\end{document}